\documentclass[12pt,oneside,reqno]{amsart}

\usepackage{fullpage}
\RequirePackage[OT1]{fontenc}
\RequirePackage{amsthm,amsmath,natbib}
\RequirePackage[colorlinks,citecolor=blue,urlcolor=blue]{hyperref}

\newcommand{\mb}{\mathbb} 
\newcommand{\mc}{\mathcal} 
\newcommand{\mbf}{\mathbf}
\newcommand{\RR}{\mb{R}}
\newcommand{\PP}{\mb{P}}
\newcommand{\ZZ}{\mathbb{Z}}
\newcommand{\NN}{\mathbb{N}}
\newcommand{\EE}{\mathbb{E}}
\newcommand{\lfled}{\longrightarrow}

\newcommand{\Cov}{\operatorname{Cov}}
\newcommand{\im}{\mathrm{i}}

\newcommand{\e}{\exp}
\newcommand{\h}{\widehat}
\newcommand{\1}{\mb{I}}

\newcommand{\vn}{\mbf{n}}

\newtheorem{thm}{Theorem}
\newtheorem{prop}{Proposition}
\newtheorem{lem}[thm]{Lemma}
\newtheorem{rmk}{Remark}[section]
\begin{document}
	
\title{A dependent Lindeberg central limit theorem for cluster functionals on stationary random fields}

\author{Jos\'e G. G\'omez-Garc\'ia$^{*}$}
\address{* Normandie Universit\'e, UNICAEN, CNRS, LMNO, France.}
%\curraddr{Laboratoire de Math\'ematiques Nicolas Oresme, CNRS UMR 6139,
%	Universit\'e de Caen-Normandie BP 5186, 14032 Caen Cedex,	France}
\email{jose3g@gmail.com; jose-gregorio.gomez-garcia@unicaen.fr}
\thanks{The author was funded by the Normandy Region RIN program.}


%    Information for 3 author

%    General info

\date{Submitted on \today}


\begin{abstract}
	In this paper, we provide a central limit theorem for the finite-dimensional marginal distributions of empirical processes $(Z_n(f))_{f\in\mc{F}}$ whose index set $\mc{F}$ is a family of cluster functionals valued on blocks of values of a stationary random field. The practicality and applicability of the result depends mainly on the usual Lindeberg condition and a sequence $T_n$ which summarizes the dependence between the blocks of the random field values. Finally, as application, we use the previous result in order to show the Gaussian asymptotic behavior of the iso-extremogram estimator introduced in this paper.
\end{abstract}

\maketitle

\section{Introduction}
Recent developments in massive data processing lead us to think in a different way about certain problems in Statistics. In particular, it is of interest to develop the construction of statistics as functions of data blocs and to study their inference. On the other hand, very often, in some applications ({\it e.g.}, in extremes \cite{Davis2008} and in astronomy \cite{Long2018}) only very little data is {\it relevant} for the estimates, without forgetting that this is also hidden among a large mass of ``raw data". This brings us to the idea of thinking about clusters of data deemed {\it ``relevant"} (or type extremal, in the context of extreme value theory), where we say that two relevant values belong to two different clusters if they belong to two different blocks. Moreover, these relevant values are in the cores of the blocks, where the {\it core of a block} $B$ is defined as the smaller sub-block $\mc{C}(B)$ of $B$ that contains all the relevant values of $B$, if they exist. 

\

In the context of this work, we consider functionals which act on these clusters of relevant values and we develop useful lemmas in order to simplify the essential step to establish a Lindeberg central limit theorem for these {\it ``cluster functionals"} on stationary random fields, inspired by the works of \cite{Bardet2007}, \cite{Drees&Rootzen2010} and \cite{Gomez2018}. 

\

Precisely, let $d\in\NN$ and let us denote $\mbf{n}:=(n_1,\ldots, n_d)$, $\mbf{1}:=(1,\ldots, 1)\in\NN^d$ and $[j]:=[1 : j]$, where $[i : j]:=\{i,i+1, \ldots, j\}\subset\ZZ$. Let $X=\left\{X_{\mbf{t}}\ : \ \mbf{t} \in \NN^d\right\}$ be a $\RR^k-$valued stationary random field and let  $\mb{X}=\left\{X_{\mbf{n},\mbf{t}}\ : \ \mbf{t}\in\right.$ $\left.[n_1]\times \cdots \times [n_d]\right\}_{\mbf{n}\in\NN^d}$ be the corresponding normalized random observations from the random field $X$, defined by $X_{\vn,\mbf{t}}=L_\vn(X_\mbf{t})\1_A(X_\mbf{t})$ for some measurable functions $L_\vn : \RR^k \lfled \RR^k$, such that 
\begin{align}\label{a1}
\PP\left(\left. X_{\mbf{n}, \mbf{1}} \in \ \cdot \ \right| X_{\mbf{n}, \mbf{1}} \in A \right) \underset{\mbf{n}\to \infty}{\lfled}G(\cdot),
\end{align}
where $G$ is a non-degenerate distribution and $A\subseteq \RR^k\setminus \{0\}$ is the {\it relevance set}. Here, $\1_A(\cdot)$ denotes the usual indicator function of a subset $A$ and the tendency $\mbf{n}\to\infty$ means that $n_i\to\infty$ for all $i\in [d]$. In particular, the convergence (\ref{a1}) is fulfilled if the random vector $X_\mbf{1}$ is regularly varying. For details about regularly varying vectors one can refer to Resnick \cite{Resnick1986,Resnick1987}.

\ 

For each $i\in[d]$, let $r_i:=r_{n_i}$ be a integer value such that $r_i=o(n_i)$ and  $m_i:=\lceil n_i/r_i \rceil:=\max\left\{k\in\NN: k \leq n_i/r_i \right\}$. We define the \textbf{$d-$blocks} (or simply {\it blocks}) of $\mb{X}$ by
\begin{equation}\label{a2}
Y_{\mbf{n},j_1\ldots j_d}:=\left( X_{\mbf{n}, \mbf{t}}\right)_{\mbf{t}\in \prod_{i=1}^d [(j_i-1)r_i+1\ : \ j_i r_i]}, 
\end{equation}
where $(j_1,\ldots, j_d)\in D_{\mbf{n},d}:=\prod_{i=1}^d[m_i]$. We have thus $m_1\cdots m_d$ complete blocks $Y_{\mbf{n},j_1\ldots j_d}$, and no more than $m_1+m_2+\cdots + m_d-d+1$ incomplete ones which we will ignore. Besides, as usual, $\prod_{i=1}^{d} A_i$ denotes the Cartesian product $A_1\times \cdots \times A_d$ and, by stationarity, we will denote $Y_\mbf{n} \overset{\mc{D}}{=}Y_{\mbf{n},\mbf{1}}$ as a generic block of $\mb{X}$.

\

We are now going to formally define the core of a block, cluster functional and the empirical process of cluster functionals, which are generalizations of the definitions of \cite{Drees&Rootzen2010} to $d-$blocks. 

\

Let $y=(x_\mbf{t})_{\mbf{t}\in\prod_{i=1}^d[r_i]}$ be a $d-$block. The \textbf{core} of the block $y$ (w.r.t. the relevance set $A$) is defined as 
\[\mc{C}(y)=\left\{\begin{array}{cc}
(x_\mbf{t})_{\mbf{t}\in \prod_{i=1}^d [r_{i,I}\ : \  r_{i,S}]}, & \mbox{ if  $x_\mbf{t}\in A$  for some $\mbf{t}\in\prod_{i=1}^d [r_i]$};  \\
0,  & \mbox{otherwise},
\end{array}\right.\]
where, for each $i\in [d]$, $r_{i,I}$  and $r_{i,S}$ are defined as
\begin{align}
r_{i,I}&=\min\left\{j_i\in [r_i] :  \ x_{(j_1,\ldots, j_i,\ldots, j_d)} \in A , \text{for some } (j_1,\ldots,j_{i-1},j_{i+1},\ldots, j_d)\in\hspace{-0.4cm} \prod_{k\in [d]\setminus \{i\}}[r_k]\right\},\nonumber
\\
r_{i,S}&=\max\left\{j_i\in [r_i] : \ x_{(j_1,\ldots, j_i,\ldots, j_d)} \in A , \text{for some } (j_1,\ldots,j_{i-1},j_{i+1},\ldots, j_d)\in\hspace{-0.4cm} \prod_{k\in [d]\setminus \{i\}}[r_k]\right\}.\nonumber
\end{align}
Let $(E, \mc{E})$ be a measurable subspace of $(\RR^k, \mc{B}(\RR^k))$ for some $k\geq1$ such that $0\in E$. Let $\mb{B}_{l_1, \ldots, l_d}(E)$ be the set of $E-$valued blocks (or arrays) of size $l_1 \times l_2 \times\cdots\times l_d$, with $l_1,\ldots, l_d \in \NN$. Consider now the set 
$$E_\cup:= \bigcup_{l_1,\ldots,l_d=1}^\infty \mb{B}_{l_1, \ldots, l_d}(E),$$
which is equipped with the $\sigma-$field $\mc{E}_\cup$ induced by the Borel$-\sigma-$fields on $\mb{B}_{l_1, \ldots, l_d}(E)$, for $l_1,\ldots, l_d \in\NN$. A \textbf{cluster functional} %$^($\footnote{This definition is only a generalization of cluster functionals of time series, with respect to the $d-$blocks and their cores.}$^)$ 
is a measurable map $f: (E_\cup, \mc{E}_\cup)\lfled (\RR, \mc{B}(\RR))$ such that 
\begin{align}\label{fc}
f(y)=f(\mc{C}(y)), \mbox{ for all }y\in E_\cup,\quad \mbox{ and } \quad \ f(0)=0.
\end{align}

\

Let $\mc{F}$ be a class of  cluster functionals and let $\left\{Y_{\mbf{n},j_1j_2\ldots j_d} \ : \ (j_1,\ldots,j_d)\in D_{\mbf{n},d}\right\}$ be the family of blocks of size $r_1\times r_2\times \cdots \times r_d$ defined in (\ref{a2}). The \textbf{empirical process $Z_\mbf{n}$ of cluster functionals} in $\mc{F}$, is the process $(Z_\mbf{n}(f))_{f\in\mc{F}}$ defined by
\begin{align}\label{PE}
Z_\mbf{n}(f):= \frac{1}{\sqrt{n_\mbf{n} v_\mbf{n}}}\sum_{(j_1,\ldots ,j_d)\in D_{\mbf{n},d}}(f(Y_{\mbf{n},j_1\ldots j_d})-\EE f(Y_{\mbf{n},j_1\ldots j_d})),
\end{align}
where  $n_\mbf{n}=n_1 \cdots n_d$ and $v_\mbf{n}:=\PP(X_{\mbf{n},\mbf{1}} \in A)$ with $A\subseteq E\setminus \{0\}$ denoting the relevance set.
%{\color{red}\huge Examples}

\

Under the Lindeberg condition and the convergence to zero of a sequence $T_\mbf{n}$ that summarizes the dependence between the blocks of values of the random field, we prove that the finite-dimensional marginal distributions (\textbf{fidis}) of the empirical process (\ref{PE}) converge to a Gaussian process. The proof basically consists of the ``Lindeberg method" as in \cite{Bardet2007}, but adapted here to stationary random fields. 

\

Regarding the condition $T_\mbf{n}\lfled 0$, as $\mbf{n}\to\infty$, this can be fulfilled if the random field $X$ has short range dependence properties, {\it e.g.}, if the random field $X$ is weakly dependent in the sense of Doukhan \& Louhichi \cite{Doukhan1999} under convenient conditions for the decay rates of the weak-dependence coefficients. These rates are calculated in \cite{Gomez2018} in the context of extreme clusters of time series.

\

The rest of the paper consists of two sections. In Section \ref{Res}, we provide useful lemmas in order to establish the central limit theorem for the fidis of the cluster functionals empirical process (\ref{PE}). In Section \ref{appli} we introduce the iso-extremogram ({\it a correlogram for extreme values of space-time processes}) and we use the CLT of Section \ref{Res} in order to show that, under additional suitable conditions, the iso-extremogram estimator has asymptotically a Gaussian behavior.
%%%%%%%%%%%%%%%%%%%%%%%%%%%%%%%%%%%%%%%%%%%%%%%%%%%%%%%%%%%%%%%%%%%%%%%%%%%%%
\section{Results}\label{Res}
In this section we provide useful lemmas that simplify notably the essential step to establish a central limit theorem for the fidis of the empirical process defined in (\ref{PE}). The proof consists in the same techniques that Bardet {\it et al. }\cite{Bardet2007} used in the demonstrations of their dependent and independent Lindeberg lemmas, but generalized here to random fields. 

\

In order to establish the CLT, firstly consider the following basic assumption:
\begin{itemize}
	\item[\textbf{(Bas)}] The vector $\mbf{r}=(r_1, \ldots, r_d)\in \NN^d$ is such that  $r_i  \ll n_i$ for each $i \in [d]$.
	\\
	Besides, denoting $r_\mbf{n}=r_1\cdots r_d$, $r_\mbf{n} v_\mbf{n} \lfled \tau < \infty$ and $n_\mbf{n} v_\mbf{n} \lfled \infty$, as $\mbf{n}\to \infty$. 
\end{itemize}
Secondly, consider the following essential convergence assumptions:
\begin{itemize}
	\item[\textbf{(Lin)}] $(r_\mbf{n}v_\mbf{n})^{-1}\EE\left[\left(f(Y_{\mbf{n}})-\EE f(Y_{\mbf{n}})\right)^2\1_{
		\left\{
		|f(Y_{\mbf{n}})-\EE f(Y_{\mbf{n}})|>\epsilon\sqrt{n_\mbf{n}v_\mbf{n}}
		\right\}
	}
	\right]=o(1)$, \ $\forall\epsilon>0$, \ $\forall f\in\mc{F}$;
	\item[\textbf{(Cov)}] $(r_\mbf{n}v_\mbf{n})^{-1}
	\Cov\left(f(Y_{\mbf{n}}),g(Y_{\mbf{n}})\right)\lfled c(f,g)$, \ $\forall f,g\in \mc{F}$.
\end{itemize}

\

Consider now the random blocks $Y_{\mbf{n},j_1 \ldots j_d}$, with $(j_1, \ldots, j_d)\in D_{\mbf{n},d}$ defined in (\ref{a2}). For each $k-$tuple of cluster functionals $\mbf{f}_k=(f_1, \ldots, f_k)$ and each $(j_1, \ldots, j_d)\in D_{\mbf{n},d}$, we define the random vector:
\begin{equation}\label{W}
%\hspace{-0.279cm}
W_{\mbf{n},j_1\ldots j_d}:=\frac{1}{\sqrt{n_\mbf{n}v_\mbf{n}}}\left(f_1(Y_{\mbf{n},j_1\ldots j_d})-\EE f_1(Y_{\mbf{n}, j_1\ldots j_d}), \ldots, f_k(Y_{\mbf{n},j_1\ldots j_d})-\EE f_k(Y_{\mbf{n}, j_1\ldots j_d})  \right).
\end{equation}

Without loss of generality and in order to simplify writing, we will consider $d=2$ in the rest of this section. 

\

Let $(W'_{\vn,ij})_{(i,j)\in D_{\vn,2}}$ be a sequence of zero mean independent $\RR^k$-valued random variables, independents of the sequence $(W_{\vn, ij})_{(i,j)\in D_{\vn, 2}}$, such that $W'_{\vn, ij} \sim \mc{N}_k\left(\mbf{0}, \Cov(W_{\vn, ij}) \right)$, for all $(i,j)\in D_{\vn, 2}$. Denote by $\mc{C}_b^3$ the set of bounded functions $h:\RR^k \lfled \RR$ with bounded and continuous partial derivatives up to order $3$. For $h\in\mc{C}_b^3$ and $\mbf{n}=(n_1, n_2)\in\NN^2$, define
\begin{equation}\label{delta}
\Delta_\mbf{n}:= \left|\EE\left[h\left(\sum_{(i,j)\in D_{\vn, 2}} W_{\vn, ij}\right) - h\left(\sum_{(i,j)\in D_{\vn, 2}} W'_{\vn, ij}\right)\right] \right|.
\end{equation}

\

The following assumption will allow us to present, in a useful and simplified form, lemmas of Lindeberg under independence and dependence.
\begin{itemize}
	\item[\textbf{(Lin')}] It exists $\delta \in (0,1]$ such that, for all $(i,j)\in D_{\vn, 2}$, $\EE\left\| W_{\vn, ij}\right\|^{2+\delta}<\infty$ for all $\vn\in \NN^2$ and all $k-$tuple of cluster functionals $(f_1, \ldots, f_k)\in\mc{F}^k$. Moreover, denote 
	$$A_\mbf{n}:= \sum_{(i,j)\in D_{\vn,2}} \EE\left\|W_{\vn,ij}\right\|^{2+\delta}.$$
\end{itemize}

\begin{lem}[Lindeberg under independence]
	\label{lem1}
	Suppose that the blocks $(Y_{\vn, ij})_{(i,j)\in D_{\vn, 2}}$ are independents and that the random variables $(W_{\vn, ij})_{(i,j)\in D_{\vn, 2}}$ defined in (\ref{W}) satisfy Assumption (Lin'). Then, for all $\vn\in\NN^2$: 
	$$\Delta_\vn\leq 6\ \| h^{(2)}\|_\infty^{1-\delta} \ \| h^{(3)}\|_\infty^\delta \ A_\vn.$$
\end{lem}
\textbf{Proof.}
	First, notice that
	\begin{equation}\label{descomp_lem1_c5} 
	\Delta_\vn \leq \sum_{(i,j)\in D_{\vn,2}} \Delta_{\vn, ij} \ ,
	\end{equation}
	where 
	%\vspace{-0.6cm}
	\begin{align}
	\Delta_{\vn, ij}&:= \left|\EE\left[h_{ij}(V_{\vn,ij}+ W_{\vn, ij})-h_{ij}(V_{\vn, ij}+W'_{\vn,ij})\right]\right|, \quad \forall (i,j)\in D_{\vn, 2} \ ;\nonumber
	\\
	V_{\vn,ij}&:=\sum_{(u,v)\in D_{\vn,2}\setminus\left(\bigcup_{l=0}^{i-1}L_l^{m_2}\cup L_i^j\right)}W_{\vn,uv}, \quad \forall (i,j)\in D_{\vn, 2}\setminus \{(m_1,m_2)\} \ , \nonumber
	\\V_{\vn, m_1m_2}&=0 \ ;\text{ and }\nonumber
	\\
	h_{ij}(x)&:=\EE \left[ h\left(x+ \sum_{u=0}^{i-1}\sum_{v=1}^{m_2} W'_{\vn,uv} + \sum_{v=0}^{j-1} W'_{\vn, iv}\right) \right].\nonumber
	\end{align}
	Besides, we set the convention $W_{\vn, ij}=0$, if either $i=0$ or $j=0$. 
	
	\
	
	Now, we will use some lines of the proof of Lemma 1 in \cite{Bardet2007}. 
	\\
	Let $v,w \in\RR^k$. From Taylor's formula, there exist vectors $v_{1,w}, v_{2,w}\in \RR^k$ such that: 
	\begin{align}
	h(v+w)&= h(v)+h^{(1)}(v)(w)+\frac{1}{2}h^{(2)}(v_{1,w})(w,w) \nonumber
	\\&= h(v)+h^{(1)}(v)(w)+\frac{1}{2}h^{(2)}(v)(w,w)+\frac{1}{6}h^{(3)}(v_{2,w})(w,w,w) \ , \nonumber
	\end{align} 
	where, for $j=1,2,3$, $h^{(j)}(v)(w_1, w_2, \ldots, w_j)$ stands for the value of the symmetric $j-$linear form from $h^{(j)}$ of $(w_1, \ldots, w_j)$ at $v$. Moreover, denote
	$$\|h^{(j)}(v)\|_1=\sup_{\| w_1\|, \ldots, \|w_j\|\leq 1} |h^{(j)}(v)(w_1,\ldots, w_j)| \quad \text{and} \quad \|h^{(j)} \|_\infty = \sup_{v\in\RR^k} \|h^{(j)}(v)\|_1.$$
	Thus, for $v,w, w' \in \RR^k$, there exist some suitable vectors $v_{1,w}, v_{2,w}, v_{1,w'}, v_{2,w'}\in\RR^k$ such that 
	\begin{multline*}
	h(v+w)-h(v+w')=h^{(1)}(v)(w-w')+\frac{1}{2}\left(h^{(2)}(v)(w,w)-h^{(2)}(v)(w',w')\right)
	\\+\frac{1}{2}\left(\left(h^{(2)}(v_{1,w})-h^{(2)}(v)\right)(w,w)-\left(h^{(2)}(v_{1,w'})-h^{(2)}(v)\right)(w',w') \right),
	\end{multline*}
	by using the approximation of Taylor of order $2$, and 
	\begin{multline*}
	h(v+w)-h(v+w')=h^{(1)}(v)(w-w')+ \frac{1}{2}\left(h^{(2)}(v)(w,w)-h^{(2)}(v)(w',w')\right)
	\\+\frac{1}{6}\left( h^{(3)}(v_{2,w})(w,w,w) - h^{(3)}(v_{2,w'})(w',w',w')\right),
	\end{multline*}
	by using the approximation of Taylor of order $3$.
	\\
	Thus, $\gamma=h(v+w)-h(v+w')-h^{(1)}(v)(w-w')-\frac{1}{2}\left(h^{(2)}(v)(w,w)-h^{(2)}(v)(w',w')\right)$ satisfies: 
	\begin{align}\label{gamma_c5}
	%\hspace{-0.31cm}
	|\gamma| &\leq \left(\left(\|w\|^2+\|w'\|^2\right)\|h^{(2)}\|_\infty\right) \wedge \left(\frac{1}{6}\left(\|w\|^3+\|w'\|^3\right)\|h^{(3)}\|_\infty\right)\nonumber
	\\
	&\leq \left(\|w\|^2 \|h^{(2)}\|_\infty \right)\wedge \left(\frac{1}{6}\|w\|^3\|h^{(3)}\|_\infty\right) 
	+  \left(\|w\|^2 \|h^{(2)}\|_\infty \right)\wedge \left(\frac{1}{6}\|w'\|^3\|h^{(3)}\|_\infty\right)\nonumber
	\\
	& +\left(\|w'\|^2 \|h^{(2)}\|_\infty \right)\wedge \left(\frac{1}{6}\|w\|^3\|h^{(3)}\|_\infty\right)
	+ \left(\|w'\|^2 \|h^{(2)}\|_\infty \right)\wedge \left(\frac{1}{6}\|w'\|^3\|h^{(3)}\|_\infty\right)\nonumber
	\\
	&\leq \frac{1}{6^\delta}\|h^{(2)}\|^{1-\delta}_\infty \|h^{(3)}\|^\delta_\infty \left(\|w\|^{2+\delta} + \|w\|^{2(1-\delta)} \|w'\|^{3\delta} + \|w\|^{3\delta} \|w'\|^{2(1-\delta)} + \|w'\|^{2+\delta}\right),
	\end{align}
	where (\ref{gamma_c5}) is given by using the inequality $1\wedge a \leq a^\delta$, with $a\geq 0$ and $\delta\in [0,1]$.
	
	\
	
	Substituting $h_{ij},V_{\vn, ij}, W_{\vn, ij}$ and $W'_{\vn, ij}$ for $h, v,w$ and $w'$ in the preceding inequality (\ref{gamma_c5}) and taking expectations, we will obtain a bound for $\Delta_{\vn, ij}$. Indeed, 
	\begin{multline*}
	\EE\left[h_{ij}(V_{\vn, ij} + W_{\vn, ij})-h_{ij}(V_{\vn, ij}+W'_{\vn, ij})\right]=\EE\left[h_{ij}(V_{\vn, ij} + W_{\vn, ij})-h_{ij}(V_{\vn, ij}+W'_{\vn, ij})\right]+0\nonumber
	\\
	=\EE\left[h_{ij}(V_{\vn, ij} + W_{\vn, ij})-h_{ij}(V_{\vn, ij}+W'_{\vn, ij})\right]- \EE\left[h^{(1)}_{ij}(V_{\vn,ij})(W_{\vn,ij}-W'_{\vn,ij})\right]
	\\
	-\frac{1}{2}\EE\left[h_{ij}^{(2)}(V_{\vn,ij})(W_{\vn,ij},W_{\vn,ij})-h_{ij}^{(2)}(V_{\vn, ij})(W'_{\vn,ij}, W'_{\vn,ij})\right] ,
	\end{multline*}
	because $V_{\vn,ij}$ is independent of $W_{\vn,ij}$ and $W'_{\vn, ij}$ , and because $\EE W_{\vn, ij}=\EE W'_{\vn, ij}=0$ and $\Cov(W_{\vn, ij})=\Cov(W'_{\vn, ij})$ for all $(i,j)\in D_{\vn, 2}$.
	\\	
	On the other hand, using Jensen's inequality, we derive 
	$\EE\|W'_{\vn, ij}\|^{2+\delta}\leq \left(\EE\|W'_{\vn,ij}\|^4 \right)^{\frac{1}{2}+\frac{\delta}{4}}$, and $\EE\|W'_{\vn, ij}\|^4\leq 3 \cdot \left( \EE\|W_{\vn, ij}\|^2\right)^2$ 
	because $W'_{\vn, ij}$ is a Gaussian random variable with the same covariance as $W_{\vn, ij}$.
	\\ 
	Therefore, 
	\begin{equation}\label{ine1_c5}
	\EE\|W'_{\vn,ij}\|^{2+\delta}
	\leq \left(3 \cdot \left( \EE\|W_{\vn, ij}\|^2\right)^2 \right)^{\frac{1}{2}+\frac{\delta}{4}}
	=3^{\frac{1}{2}+\frac{\delta}{4}}\left( \EE\|W_{\vn, ij}\|^2\right)^{1+\frac{\delta}{2}}
	\leq 3^{\frac{1}{2}+\frac{\delta}{4}} \EE\|W_{\vn,ij}\|^{2+\delta},
	\end{equation}
	\begin{align}\label{ine2_c5}
	\hspace{-0.28cm}
	\EE\|W'_{\vn, ij}\|^{2(1-\delta)}\EE\|W_{\vn, ij}\|^{3\delta}&\leq \left(\EE\|W'_{\vn, ij}\|^2\right)^{1-\delta}\EE\|W_{\vn,ij}\|^{3\delta}
	\nonumber
	\\
	&\leq \left(\EE\|W_{\vn, ij}\|^2\right)^{1-\delta}\EE\|W_{\vn,ij}\|^{3\delta}\leq \EE\|W_{\vn,ij}\|^{2+\delta} \ . \
	\end{align}
	Besides, for $3\delta<2$, 
	\begin{align}\label{ine3_c5}
	\EE\|W_{\vn, ij}\|^{2(1-\delta)} \EE\|W'_{\vn, ij}\|^{3\delta}\leq \EE\|W_{\vn, ij}\|^{2(1-\delta)}\left( \EE\|W'_{\vn, ij}\|^2\right)^{\frac{3\delta}{2}}\leq \EE\|W_{\vn, ij}\|^{2+\delta},
	\end{align}
	else 
	\begin{align}\label{ine4_c5}
	\hspace{-0.19cm}
	\EE\|W_{\vn, ij}\|^{2(1-\delta)}\EE\|W'_{\vn, ij}\|^{3\delta} &\leq \EE\|W_{\vn, ij}\|^{2(1-\delta)}\left( \EE\|W'_{\vn, ij}\|^4\right)^{\frac{3\delta}{4}}, \quad \text{ because }3\delta\leq 4\nonumber
	\\
	& \leq 3^{\frac{3\delta}{4}} \EE\|W_{\vn, ij}\|^{2(1-\delta)}\left(\EE\|W_{\vn, ij}\|^{2}\right)^{\frac{3\delta}{2}} \leq 3^{\frac{1}{2}+ \frac{\delta}{4}}\EE \| W_{\vn, ij}\|^{2+\delta}.
	\end{align}
	The inequalities (\ref{ine1_c5})-(\ref{ine4_c5}) allow to simplify the terms between parentheses in the last inequality in (\ref{gamma_c5}). Recall that $\|h_{ij}^{(k)}\|_\infty\leq \|h^{(k)}\|_\infty$ for all $(i,j)\in D_{\vn,2}$ and $0\leq k\leq 3$.
	Therefore, we obtain that 
	\begin{align}
	\Delta_{\vn,ij}& \leq \frac{2(1+3^{\frac{1}{2}+\frac{\delta}{4}})}{6^\delta}\|h^{(2)}\|_\infty^{1-\delta}\|h^{(3)}\|_\infty^\delta \EE\|W_{\vn, ij}\|^{2+\delta} \leq 6 \|h^{(2)}\|_\infty^{1-\delta}\|h^{(3)}\|_\infty^\delta \EE\|W_{\vn, ij}\|^{2+\delta}, \nonumber
	\end{align}
	because, for all $\delta\in [0,1]$, \ $C(\delta)=\frac{2(1+3^{\frac{1}{2}+\frac{\delta}{4}})}{6^\delta}\leq C(0)=2(1+\sqrt{3})<6$.
	\\
	\\
	As a consequence, from Assumption (Lin'), $\Delta_{\vn}\leq 6 \ \| h^{(2)}\|_\infty^{1-\delta} \ \| h^{(3)}\|_\infty^\delta \ A_\vn$. \hspace{0.5cm}%$\blacksquare$
	\hfill$\square$\vskip2mm\hfill
\begin{rmk}\label{rem_cond_C1} Taking $\epsilon < 6 \|h^{(2)}\|_\infty \ (\|h^{(3)}\|_\infty)^{-1}$ and using suitably the second inequality of (\ref{gamma_c5}) in the proof of Lemma \ref{lem1}, classical Lindeberg conditions may be used:
	\begin{align}\label{Cota_lindeberg}
	\Delta_{\vn}\leq 2 \ \|h^{(2)}\|_\infty\ B_\vn(\epsilon)+\|h^{(3)}\|_\infty \ a_\vn \ \left(\frac{4}{3}\epsilon +\sqrt{B_\vn(\epsilon)}\right),
	\end{align}
	where
	\begin{align}
	B_\vn(\epsilon)&=\sum_{(i,j)\in D_{\vn, 2}} \EE\left[\|W_{\vn,ij}\|^2\1_{\{\|W_{\vn,ij}\|>\epsilon\}}\right], \quad \epsilon>0, \  \vn\in\NN^2 \ ;\nonumber
	\\
	a_\vn&=\sum_{(i.j)\in D_{\vn,2}} \EE\|W_{\vn,ij}\|^2< \infty, \quad  \vn\in\NN^2.\nonumber
	\end{align}
	Moreover, these classical Lindeberg conditions derive the conditions from Lemma \ref{lem1}. Indeed, 
	$$\Delta_\vn \leq 2 \ \|h^{(2)}\|_\infty \ \epsilon^{-\delta} A_\vn +\|h^{(3)}\|_\infty \ a_\vn \ \left(\frac{4}{3}\epsilon +\epsilon^{-\delta/2} \sqrt{A_\vn}\right),$$
	for  $\delta\in (0,1)$ and $\epsilon> 0$.
\end{rmk}
The proof of this remark for general independent random vectors is given in \cite[p.165]{Bardet2007}.
\begin{rmk}\label{rem2}	Observe that the assumptions (Lin) and (Cov) imply that $B_\vn(\epsilon) \underset{\vn\to\infty}{\lfled}0$ and that $a_\vn=\sum_{i=1}^{k}(r_\mbf{n}v_\vn)^{-1}\Cov\left(f_i(Y_\vn), f_i(Y_\vn)\right)\underset{\vn\to\infty}{\lfled}  \sum_{i=1}^k c(f_i,f_i)<\infty$, respectively. Therefore, if the blocks $(Y_{\vn, ij})_{(i,j)\in D_{\vn, 2}}$ are independent and if the assumptions (Lin) and (Cov) hold, then from Lemma \ref{lem1} and Remark \ref{rem_cond_C1}, the fidis of the empirical process $(Z_\vn (f))_{f\in\mc{F}}$ of cluster functionals converge to the fidis of a Gaussian process $(Z(f))_{f\in\mc{F}}$ with covariance function $c$. %The proof is immediate from Lemma \ref{lem1_c5} and Remark \ref{rem_cond_C1_c5}.
\end{rmk}

For the dependent case, we need to consider more notations: 
\\
Let $L_i^j:=\left\{(i,v) \ : \ v\in[j] \right\}\subset D_{\mbf{n},2}$ , for all $(i,j)\in D_{\mbf{n},2}$. We set $L_i^0=L_0^j=\emptyset$ for any $i\in [m_1]$ and any $j\in [m_2]$. For each $k\in\NN$, \ $\mbf{f}_k=(f_1, \ldots, f_k)\in\mc{F}^k$, \ $\mbf{t}\in\RR^k$ and $\vn\in\NN^2$ ; we define
\begin{equation*}\label{T}
T_{\mbf{n},\mbf{t}}(\mbf{f}_k):=\hspace{-0.4cm}\sum_{(j_1,j_2)\in D_{\mbf{n},2}}\left|\Cov\left(\exp\left(\im \langle \mbf{t}, \hspace{-0.4cm}\sum_{(u_1,u_2)\in D_{\mbf{n},2}\setminus (\bigcup_{l=0}^{j_1-1}L_l^{m_2} \cup L_{j_1}^{j_2})}\hspace{-0.4cm}W_{\mbf{n},u_1 u_2}\rangle\right),\exp\left(\im \langle\mbf{t},W_{\mbf{n},j_1j_2}\rangle \right)\right)\right|.
\end{equation*}
\begin{lem}[Dependent Lindeberg lemma]\label{lem2} 
	Suppose that the r.v.'s $(W_{\vn, ij})_{(i,j)\in D_{\vn, 2}}$ defined in (\ref{W}) satisfy Assumption (Lin'). Consider the special case of complex exponential functions $h(\mbf{w})=\e\left(\im \langle \mbf{t}, \mbf{w}\rangle\right)$ with $\mbf{t}\in\RR^k$. Then, for each $k \in \NN$ and each $k-$tuple $\mbf{f}_k=(f_1, \ldots, f_k)$ of cluster functionals, the following inequality holds:
	$$\Delta_\vn\leq T_{\vn, \mbf{t}}(\mbf{f}_k) + 6 \| \mbf{t}\|^{2+\delta} A_\vn \ , \qquad \vn\in\NN^2.$$
	%where $T_{\vn, \mbf{t}}(\mbf{f}_k)$ is defined in (\ref{TT}).
\end{lem}
\textbf{Proof.}
	Consider $(W^*_{\vn, j_1j_2})_{(j_1,j_2)\in D_{\vn,2}}$ an array of independent random variables satisfying Assumption (Lin') and such that $(W^*_{\vn, j_1j_2})_{(j_1,j_2)\in D_{\vn,2}}$ is independent of $(W_{\vn, j_1j_2})_{(j_1,j_2)\in D_{\vn,2}}$ and $(W'_{\vn, j_1j_2})_{(j_1,j_2)\in D_{\vn,2}}$. Moreover, assume that $W^*_{\vn, j_1j_2}$ has the same distribution as $W_{\vn, j_1j_2}$ for $(j_1,j_2)\in D_{\vn,2}$. 
	\\
	Then, using the same decomposition (\ref{descomp_lem1_c5}) in the proof of the previous lemma, one can also write,
	\begin{multline}\label{descomp2_c5}
	\Delta_{\vn, j_1j_2}\leq \left|\EE\left[h_{j_1j_2}(V_{\vn,j_1j_2}+ W_{\vn, j_1j_2})-h_{j_1j_2}(V_{\vn, j_1j_2}+W^*_{\vn,j_1j_2})\right]\right|
	\\
	+\left|\EE\left[h_{j_1j_2}(V_{\vn,j_1j_2}+ W^*_{\vn, j_1j_2})-h_{j_1j_2}(V_{\vn, j_1j_2}+W'_{\vn,j_1j_2})\right]\right|.
	\end{multline}
	%\\
	%Consider now the special case of complex exponential functions $h(\mbf{w})=\exp\left(\im \langle \mbf{t}, \mbf{w}\rangle\right)$ for $\mbf{t}\in\RR^k$, where $\langle \mbf{u}, \mbf{v}\rangle$ denotes the scalar product in $\RR^k$. 
	Then, from the previous lemma, the second term of the RHS of the inequality (\ref{descomp2_c5}) is bounded by 
	$$ 6\ \|h^{(2)}\|_\infty^{1-\delta} \ \|h^{(3)}\|_\infty^\delta \ \EE\|W_{\vn, j_1j_2}\|^{2+\delta} \leq 6\ \| \mbf{t}\|^{2+\delta}\ \EE\|W_{\vn, j_1j_2}\|^{2+\delta}.$$
	For the first term of the RHS of the inequality (\ref{descomp2_c5}), first notice that for a $\RR^k-$valued random vector $X$ independent from $(W'_{\vn, j_1j_2})_{(j_1,j_2)\in D_{\vn, 2}}$, 
	\begin{align}
	\EE h_{j_1j_2}(X)&=\EE \left[ h\left(X+ \sum_{u=0}^{j_1-1}\sum_{v=1}^{m_2} W'_{\vn,uv} + \sum_{v=0}^{j_2-1} W'_{\vn, j_1v}\right) \right]\nonumber
	\\
	&=\exp\left( -\frac{1}{2} \mbf{t}^T \left(\sum_{u=0}^{j_1-1}\sum_{v=1}^{m_2} C_{\vn,uv}+\sum_{v=0}^{j_2-1}C_{\vn, j_1v} \right) \mbf{t}\right)\cdot\EE\left[\exp \left(\im \langle \mbf{t}, X\rangle \right)\right],\nonumber
	\end{align}
	because $W'_{\vn, j_1j_2} \sim \mc{N}_k\left(\mbf{0}, C_{\vn, j_1j_2} \right)$, where $C_{\vn,j_1j_2}:=\Cov(W_{\vn, j_1j_2})$ is the covariance matrix of the vector $W_{\vn, j_1j_2}$,  \ for $(j_1,j_2)\in D_{\vn, 2}$. For $j_1=0$ or $j_2=0$, recall that $W_{\vn, j_1j_2}=0$. In this case, we also set $C_{\vn,j_1j_2}=0$. 
	\\
	Thus, 
	\begin{multline*}
	\left|\EE\left[h_{j_1j_2}(V_{\vn,j_1j_2}+ W_{\vn, j_1j_2})-h_{j_1j_2}(V_{\vn, j_1j_2}+W^*_{\vn,j_1j_2})\right]\right|
	\\
	=\left|\exp\left( -\frac{1}{2} \mbf{t}^T \left(\sum_{u=0}^{j_1-1}\sum_{v=1}^{m_2} C_{\vn,uv}+\sum_{v=0}^{j_2-1}C_{\vn, j_1v} \right) \mbf{t}\right)\right.\qquad \qquad \qquad \qquad \qquad \qquad
	\\
	 \times \left.\EE\left[\exp \left(\im \langle \mbf{t}, V_{\vn, j_1j_2}\rangle \right) \cdot \left(\exp \left(\im \langle \mbf{t}, W_{\vn, j_1j_2}\rangle \right) - \exp \left(\im \langle \mbf{t}, W^*_{\vn, j_1j_2}\rangle\right)\right)\right]\right|
	\\
	=\left|\exp\left( -\frac{1}{2} \mbf{t}^T \left(\sum_{u=0}^{j_1-1}\sum_{v=1}^{m_2} C_{\vn,uv}+\sum_{v=0}^{j_2-1}C_{\vn, j_1v} \right) \mbf{t}\right)\right| \qquad \qquad \qquad \qquad \qquad \qquad 
	\\
	\times \left|\Cov\left(\exp \left(\im \langle \mbf{t}, V_{\vn, j_1j_2}\rangle \right),\exp \left(\im \langle \mbf{t}, W_{\vn, j_1j_2}\rangle \right) \right)\right|
	\\
	\leq \left|\Cov\left(\exp \left(\im \langle \mbf{t}, V_{\vn, j_1j_2}\rangle \right),\exp \left(\im \langle \mbf{t}, W_{\vn, j_1j_2}\rangle \right) \right)\right|.
	\end{multline*}
	Therefore, 
	\begin{align}
	\Delta_\vn 
	&= \sum_{(j_1, j_2)\in D_{\vn, 2}} \Delta_{\vn, j_1 j_2}
	\nonumber
	\\
	&\leq \sum_{(j_1, j_2)\in D_{\vn, 2}}\left(\left|\Cov\left(\exp \left(\im \langle \mbf{t}, V_{\vn, j_1j_2}\rangle \right),\exp \left(\im \langle \mbf{t}, W_{\vn, j_1j_2}\rangle \right) \right)\right| + 6 \| \mbf{t}\|^{2+\delta}\EE\|W_{\vn, j_1j_2}\|^{2+\delta}\right)
	\nonumber
	\\
	&= T_{\vn, \mbf{t}}(\mbf{f}_k) + 6\| \mbf{t}\|^{2+\delta} A_\vn. \hspace{0.5cm}%\blacksquare
	\nonumber
	\end{align}
\hfill$\square$\vskip2mm\hfill

The previous lemma together with Remark \ref{rem_cond_C1} derive the following theorem.
\begin{thm}[CLT for cluster functionals on random fields]\label{theo1_c5}
	Suppose that the basic assumption (Bas) holds and that the assumptions (Lin) and (Cov) are satisfied. Then, if for each $k\in\NN$, $T_{\mbf{n},\mbf{t}}(\mbf{f}_k)$ converges to zero as $\mbf{n}\to\infty$, for all $\mbf{t}\in\RR^k$ and all $k-$tuple $\mbf{f}_k=(f_1, \ldots, f_k)\in\mc{F}^k$ of cluster functionals, the fidis of the empirical process $(Z_\mbf{n}(f))_{f\in\mc{F}}$ of cluster functionals converge to the fidis of a Gaussian process $(Z(f))_{f\in\mc{F}}$ with covariance function $c$ defined in (Cov).
\end{thm}
\textbf{Proof.}
	The assumptions (Lin) and (Cov) imply that, as $\mbf{n}\to \infty$,  $B_\vn(\epsilon)\lfled 0$ and $a_\vn \lfled  \sum_{s=1}^k c(f_s,f_s)<\infty$, respectively. Therefore, taking into account Remark \ref{rem_cond_C1}, we obtain from Lemma \ref{lem2} that, for each $k\in\NN$,
	\begin{equation*}
	\Delta_\vn= \left|\EE\left[h\left(\sum_{(i,j)\in D_{\vn, 2}} W_{\vn, ij}\right) - h\left(\sum_{(i,j)\in D_{\vn, 2}} W'_{\vn, ij}\right)\right] \right|  \underset{\vn\to\infty}{\lfled}0,
	\end{equation*}
	for all $\mbf{t}\in\RR^k$, with $h(\mbf{w})=\exp(\im \langle \mbf{t},  \mbf{w}\rangle)$, because  by hypothesis, $T_{\vn,\mbf{t}} (\mbf{f}_k)\underset{\vn\to\infty}{\lfled} 0$ for all $\mbf{t}\in\RR^k$ and all $\mbf{f}_k=(f_1, \ldots, f_k)\in\mc{F}^k$.
	\\
	Notice that $$W'_\vn:=\sum_{(i,j)\in D_{\vn, 2}} W'_{\vn, ij} \sim \mc{N}_k(\mbf{0}, m_1m_2\Cov(W_{\vn, 11}))$$ and that $\left|\EE\left(h(W'_\vn)-h(W) \right) \right|\underset{\vn\to\infty}{\lfled}0$, where $W \sim \mc{N}_k(\mbf{0}, \Sigma_k)$, with $\Sigma_k=(c(f_i,f_j))_{(i,j)\in [k]^2}$.
	\\
	\\
	Using triangular inequality, we deduce that $$\left|\EE\left[h\left(\sum_{(i,j)\in D_{\vn, 2}} W_{\vn, ij}\right)-h(W) \right]\right|\underset{\vn\to\infty}{\lfled}0,$$ and therefore $(Z_\vn(f_1), \ldots, Z_\vn(f_k))=\sum_{(i,j)\in D_{\vn, 2}} W_{\vn, ij}\overset{\mc{D}}{\underset{\vn\to\infty}{\lfled}} W$.
	\hspace{0.5cm}%$\blacksquare$ 
\hfill$\square$\vskip2mm\hfill
\begin{rmk}\label{rem2_c5}
	The previous theorem can be formulated for $d=3$ as follows. Define $S_i=\{(u,v,w) \ : \ u \in [i],  \ v\in [m_2], \ w\in [m_3] \}\subseteq D_{\vn, 3}$, for $i\in [m_1]$, with the convention $S_{0}=\emptyset$. Moreover, $L_{ij}^k=\{(i,j,w) \ : \ w\in [k] \}$, for $(i,j,k)\in D_{\vn,3}$, and $L_{ij}^k=\emptyset$ if $i$, $j$ or $k$ is zero. Then, if (Bas), (Lin), (Cov) are satisfied (for $d=3$), and if for each $k\in\NN$,
	\begin{align}\label{T2_c5}
	&T^*_{\mbf{n},\mbf{t}}(\mbf{f}_k)=\sum_{(j_1,j_2,j_3)\in D_{\mbf{n},3}}\left|\Cov\left(\exp(\im \langle \mbf{t}, V_{\vn,j_1j_2j_3} \rangle ),\exp\left(\im \langle\mbf{t},W_{\mbf{n},j_1j_2j_3}\rangle \right)\right)\right|
	\end{align}
	converges to zero as $\mbf{n}\to\infty$ for all $\mbf{t}\in\RR^k$ and all $k-$tuple $\mbf{f}_k=(f_1, \ldots, f_k)\in\mc{F}^k$ of cluster functionals, with
	$$V_{\vn,j_1j_2j_3}:=\sum_{(u_1, u_2, u_3)\in D_{\vn, 3}\setminus \left(S_{j_1-1} \  \cup \  \bigcup_{l=0}^{j_2-1}L_{j_1 l}^{m_3}\cup L_{j_1 j_2}^{j_3} \right)} W_{\vn, u_1u_2u_3} \ ,$$
	the fidis of the empirical process $(Z_\mbf{n}(f))_{f\in\mc{F}}$ of cluster functionals converge to the fidis of a Gaussian process $(Z(f))_{f\in\mc{F}}$ with covariance function $c$. 
\end{rmk}
\begin{rmk}\label{rem1}
	We have mentioned previously that $\mbf{n}=(n_1, \ldots, n_d)\to\infty$ means $n_i\to\infty$ for each $i\in [d]$. However, if the reader would like it, the limits of the sequences indexed with $\mbf{n}$, as $\vn\to\infty$, could be reformulated in terms of the limits of such sequences as {\it ``$\mbf{n}\to\infty$ along a monotone path on the lattice $\NN^d$"}, {\it i.e.} along $\mbf{n}=(\lceil \vartheta_1(n)\rceil , \ldots, \lceil \vartheta_{d}(n)\rceil)$ for some strictly increasing continuous functions $\vartheta_i : [1,\infty) \lfled [1,\infty)$,  with $i\in [d]$, such that $\vartheta_i(n) \lfled \infty$ as $n\to\infty$, for $i\in[d]$.
\end{rmk}

Suppose that from each block $Y_\mbf{n}$ we extract a sub-block $Y'_\mbf{n}$ and that the remaining parts $R_\mbf{n}=Y_\mbf{n}-Y'_\mbf{n}$ of the blocks $Y_\mbf{n}$ do not influence the process $Z_\mbf{n}(f)$. In particular, this last statement is fulfilled if, $(r_\mbf{n} v_\mbf{n})^{-1}\EE|\Delta_\mbf{n}(f)-\EE\Delta_\mbf{n}(f)|^2\1_{\{|\Delta_\mbf{n}(f)-\EE\Delta_\mbf{n}(f)|\leq \sqrt{n_\mbf{n} v_\mbf{n}}\}}=o(1)$ and $\PP\left(|\Delta_\mbf{n}(f)-\EE\Delta_\mbf{n}(f)|>\sqrt{n_\mbf{n} v_\mbf{n}}\right)=o(r_\mbf{n}/n_\mbf{n})$, where $\Delta_\mbf{n}(f):=f(Y_\mbf{n})-f(Y'_\mbf{n})$. This assumption would allow us to consider $T_{\mbf{n}, \mbf{t}}(\mbf{f}_k)$ (or $T^*_{\mbf{n}, \mbf{t}}(\mbf{f}_k)$) as a function of the blocks $Y'_\mbf{n}$ (separated by $l_\mbf{n}$) instead of the blocks $Y_\mbf{n}$, in order to provide them bounds based on either the strong mixing coefficient of \cite{Rosenblatt1956} or the weak-dependence coefficients of \cite{Doukhan1999} for stationary random fields. These bounds are developed in \cite{Gomez2018} for the case of weakly-dependent time series, however, we will not develop them in the random field context as this is not the aim of this work. This topic will be addressed in a forthcoming applied statistics paper with numerical simulations.	
%%%%%%%%%%%%%%%%%%%%%%%%%%%%%%%%%%%%%%%%%%%%%%%%%%%%%%%%%%%%%%%%%%%%%%%%%%%%%%
\section{Asymptotic behavior of the extremogram for space-time processes}\label{appli}
In this section we propose a measure (in two versions) of serial dependence on space and time of extreme values of space-time processes. We provide an estimator for this measure and we use Theorem \ref{theo1_c5} in order to establish an asymptotic result. This section is inspired of the extremogram for times series defined in \cite{Davis2009}. 

\

Let $X=\{X_t(\mbf{s}) \ : \mbf{s}\in \ZZ^d, \ t\geq0\}$ be a $\RR^k-$valued space-time process, which is stationary in both space and time. We define the \textbf{extremogram} of $X$ for two sets $A$ and $B$ both bounded away from zero by
\begin{align}\label{extremogram_c5}
\rho_{A,B}(\mbf{s}, h_t):=\lim_{x\to\infty} \PP\left(x^{-1}X_{h_t}(\mbf{s})\in B \left| x^{-1}X_0(\mbf{0})\in A\right. \right),  \quad (\mbf{s},h_t)\in \ZZ^d\times [0,\infty) \ ;
\end{align}
provided that the limit exists. 

\

In estimating the extremogram, the limit on $x$ in (\ref{extremogram_c5}) is replaced by a high quantile $u_n$ of the process. Defining $u_n$ as the $(1-1/k_n)-$quantile of the stationary distribution of $\|X_t(\mbf{s})\|$ or related quantity, with $k_n=o(n)\lfled\infty$, as $n\to\infty$, one can redefine (\ref{extremogram_c5}) by
\begin{align}\label{extremogram2_c5}
\rho_{A,B}(\mbf{s}, h_t)=\lim_{n \to\infty} \PP\left(u_n^{-1}X_{h_t}(\mbf{s})\in B \left| u_n^{-1}X_0(\mbf{0})\in A\right. \right),  \quad (\mbf{s},h_t)\in \ZZ^d\times [0,\infty).
\end{align}
The choice of such a sequence of quantiles $(u_n)_{n\in \NN}$ is not arbitrary. The main condition to guarantee the existence of the limit (\ref{extremogram2_c5}) for any two sets $A$ and $B$ bounded away from zero,  is that it must satisfy the following convergence
\begin{align}\label{RV_c5}
k_n \PP \left( u_n^{-1} \left(X_{t_1}(\mbf{s}_1), \ldots, X_{t_p}(\mbf{s}_p) \right) \in \ \cdot \ \right) \underset{n\to\infty}{\overset{vague}{\lfled}} m_{(\mbf{s}_1,t_1),\ldots, (\mbf{s}_p, t_p)}(\ \cdot \ ),
\end{align}
for all $(\mbf{s}_i, t_i)\in \ZZ^d \times [0,\infty)$, $i\in [p]$, $p\in\NN$, where 
$$\left(m_{(\mbf{s}_1,t_1),\ldots, (\mbf{s}_p, t_p)}\right)_{(\mbf{s_i}, t_i)\in \ZZ^d\times [0,\infty), \ i \in [p], \ p\in\NN}$$
is a collection of Radon measures on the Borel $\sigma-$field $\mc{B}(\overline{\RR}^{kp}\setminus\{\mbf{0}\})$, not all of them being the null measure, with $m_{(\mbf{s}_1,t_1),\ldots, (\mbf{s}_p, t_p)}(\overline{\RR}^{kp}\setminus \RR^{kp})=0$. 
In this case, 
\begin{align}%\label{extremogram_equi}
\PP\left(u_n^{-1}X_{h_t}(\mbf{s})\in B \left| u_n^{-1}X_0(\mbf{0})\in A\right. \right)=& \dfrac{k_n\PP\left(u_n^{-1}\left(X_0(\mbf{0}), X_{h_t}(\mbf{s})\right) \in A \times B\right)}{k_n\PP\left(u_n^{-1}X_0(\mbf{0})\in A \right)}\nonumber
\\
\lfled& \dfrac{m_{(\mbf{0},0), (\mbf{s},h_t)}(A\times B)}{m_{(\mbf{0},0)}(A)} = \rho_{A,B}(\mbf{s},h_t),\nonumber 
\end{align}
provided that $m_{(\mbf{0},0)}(A)>0$.
\begin{rmk}\label{rem1_ext_c5}
	The condition (\ref{RV_c5}) is particularly satisfied if the space-time process $X$ is regularly varying. For details and examples of regularly varying space-time processes and time series, see \cite{Davis2008} and \cite{Segers2009}, respectively.  
\end{rmk}
Note that the extremogram (\ref{extremogram2_c5}) is a function of two lags: a spatial-lag $\mbf{s}\in\ZZ^d$ and a non-negative time-lag $h_t$. Due to all the spatial values that the spatial-lag $\mbf{s}$ takes, in practice, it is very complicated to analyze the results of the estimation of such extremogram. Moreover, the calculation would be computationally very slow. In order to obtain a simpler interpretation and simplify the calculations, we will assume that the space-time process $X$ satisfies the following ``isotropy" condition:
\begin{itemize}
	\item[\textbf{(I)}] For each pair of non-negative integers $h_t$ and $h_s$,
	$$\PP\left(X_0(\mbf{0})\in A, X_{h_t}(\mbf{s})\in B \right)=\PP\left(X_0(\mbf{0})\in A, X_{h_t}(\mbf{s}')\in B \right),\quad  \forall \mbf{s}, \mbf{s}' \in \mb{S}^{d-1}_{h_s} \ ,$$
\end{itemize}
where $\mb{S}^{d-1}_h:=\left\{ \mbf{s} \in\ZZ^d :  \|\mbf{s}\|_\infty=h \right\}$ with $h\geq0$ and $\|(s_1, \ldots, s_d)\|_\infty=\max_{i=1,\ldots, d} |s_i|$.
\\
Under this condition, the extremogram (\ref{extremogram2_c5})  can be redefined using only two non-negative integer lags: a spatial-lag $h_s$ and a time-lag $h_t$. Indeed, under the assumption (I), we define the \textbf{iso-extremogram} of $X$  for two sets $A$ and $B$ both bounded away from zero by
\begin{align}\label{extremogram_iso_c5}
\rho^*_{A,B}(h_s,h_t)=\rho_{A,B}(h_s \vec{e}_1, \ h_t), \qquad h_s, h_t =0,1,2,\ldots 
\end{align}
where $\vec{e}_1=(1,0,0,0,\ldots,0)\in \RR^d$ is the first element of the canonical basis of $\RR^d$.

\

We will propose now a estimator for the iso-extremogram. For this, consider w.l.o.g. $d=2$, because the case  $d>2$ is similar. 
\\
Let $X_\vn:=\left\{ X_t(i,j) : (i,j,t)\in [n_1]\times [n_2]\times [n_3]\right\}$ be the observations from a $\RR^k-$valued space-time process $X$, stationary in both space and time, and which satisfies the condition (I).  Let us set $n=n_1 n_2 n_3$. The sample iso-extremogram based on the observations $X_\vn$ is given by
\begin{multline}\label{est_extremogram_c5}
\h{\rho^*}_{A,B}(h_s, h_t):=
\dfrac{\displaystyle\sum_{(j_1,j_2)\in[m_1]\times [m_2]}\sum_{t=1}^{n_3-h_t} \sum_{(i_1,i_2)\in \mb{S}_{h_s}(c_{j_1j_2})}  \frac{\1_{\left\{\frac{X_{t+h_t}(i_1,i_2)}{u_n}\in B, \  \frac{X_t(c_{j_1j_2})}{u_n} \in A\right\}}}{\# \mb{S}_{h_s}(c_{j_1j_2})} }{\displaystyle\sum_{(j_1,j_2)\in [m_1]\times [m_2]}\sum_{t=1}^{n_3} \1_{\left\{\frac{X_t(c_{j_1j_2})}{u_n} \in A \right\}}},
\end{multline}
for $h_s=0,1,2,\ldots, \lceil 2^{-1}\min\{r_1,r_2\}\rceil-1$, and $h_t=0, \ldots, n-1$, where
$$c_{ij}:=\left(\left\lceil \frac{(2i-1)r_1+1}{2}\right\rceil, \left\lceil \frac{(2j-1)r_2+1}{2}\right\rceil\right)$$
denotes the ``center" of the block $B_{ij}=[(i-1)r_1+1: ir_1]\times[(j-1)r_2+1: jr_2]$, for $(i,j)\in [m_1]\times [m_2]$. Besides, $\mb{S}_h(u,v):=\{(i,j)\in [n_1]\times [n_2]: \|(u,v)-(i,j)\|_\infty=h\}$ with $h\geq0$ and \ $\#E$ denotes the cardinality of the set $E$. Remember that $r_i=r_{n_i}$ and $m_i=\lceil n_i/r_i\rceil$, for $i=1,2,3$.

\

Defining the cluster functional $f_{A,B,h_1,h_2}: \left(\bigcup_{l_1, l_2, l_3=1}^\infty \mb{B}_{l_1l_2l_3}(\RR^k), \mc{R}_\cup\right) \lfled (\RR,\mc{B}(\RR))$, for $h_1,h_2=0,1,2, \ldots$, such that 
\begin{equation}\label{functional_ext_c5}
\hspace{-0.5cm}
f_{A,B,h_1,h_2}\left((x_{(i_1,i_2,i_3)})_{(i_1,i_2,i_3)\in[l_1]\times [l_2]\times [l_3]}\right)=\hspace{-0.2cm}\sum_{(i_1,i_2)\in \mb{S}_{h_1}(c)}\sum_{i_3=1}^{l_3-h_2} \frac{\1_{A\times B}(x_{(c,i_3)},x_{(i_1,i_2,i_3+h_2)})}{\# \mb{S}_{h_1}(c)},
\end{equation}
with $c=(\lceil (l_1+1)/2\rceil, \lceil (l_2+1)/2\rceil)\in[l_1]\times [l_2]$ (the ``center" of the block $B=[l_1]\times[l_2]$), we can rewrite the estimator (\ref{est_extremogram_c5}) as:
\begin{align}\label{est_extremogram_funct_c5}
\h{\rho^*}_{A,B}(h_s, h_t)=\dfrac{\sum_{(j_1,j_2,j_3)\in D_{\vn,3}}f_{A,B,h_s,h_t}(Y_{\vn,j_1j_2j_3})+\delta_\vn+R_{A,B,h_s,h_t} }{\sum_{(j_1,j_2,j_3)\in D_{\vn,3}}f_{A,A,0,0}(Y_{\vn,j_1j_2j_3}) + R_{A,A,0,0}},
\end{align}
where 
\begin{align}
\delta_\vn:&=\sum_{(j_1,j_2,j_3)\in D_{\vn,3}}\sum_{(i_1,i_2)\in \mb{S}_{h_s}(c_{j_1j_2})}\sum_{t=j_3 r_3 -h_t+1}^{j_3r_3}\frac{\1_{\left\{\frac{X_{t+h_t}(i_1,i_2)}{u_n}\in B, \  \frac{X_t(c_{j_1j_2})}{u_n} \in A\right\}}}{\# \mb{S}_{h_s}(c_{j_1j_2})}, \nonumber
\\
R_{A,B,h_s,h_t}:&=\sum_{(j_1,j_2)\in [m_1]\times[m_2]}\sum_{(i_1,i_2)\in \mb{S}_{h_2}(c_{j_1j_2})}\sum_{t=m_3 r_3 + 1}^{n_3-h_t}\frac{\1_{\left\{\frac{X_{t+h_t}(i_1,i_2)}{u_n}\in B, \  \frac{X_t(c_{j_1j_2})}{u_n} \in A\right\}}}{\# \mb{S}_{h_s}(c_{j_1j_2})}\nonumber.
\end{align}
We can therefore write (\ref{est_extremogram_funct_c5}) in terms of empirical processes of cluster functionals (\ref{PE}) and use Lindeberg CLT for cluster functionals on random fields (Theorem \ref{theo1_c5}) together with suitable conditions of joint distributions, in order to prove the convergence in distribution of the iso-extremogram estimator. 

\

For this, firstly we set some considerations: the normalized random variables are defined here by $X_{\vn,(i_1,i_2,t)}=u_n^{-1}X_t(i_1,i_2)$, where $\vn=(n_1,n_2,n_3)$ and $n=n_1n_2n_3$;  and the random blocks $(Y_{\vn,j_1j_2j_3})_{(j_1,j_2,j_3)\in D_{\vn,3}}$ as in (\ref{a2}).  We define $\NN_0:=\NN\cup \{0\}$ and $\mc{F}_{A,B}:=\left\{f_{A,B,h_s,h_t} : h_s,h_t\in \NN_0 \right\}$ as the family of cluster functionals defined in (\ref{functional_ext_c5}). Moreover, for the set $A$, bounded away from zero, let $v_n:=\PP\left(u_n^{-1} X_0{(0,0)} \in A \right)$.
\\
\\
Secondly, consider the following conditions:
\begin{itemize}
	\item[\textbf{(Cov')}] For each $h_s,h_s',h_t, h_t' \in \NN_0$,
\end{itemize}
\begin{equation*}\label{cov_ext1_c5}
%\hspace{-1.3cm}
\sum_{\mbf{i}\in \mb{S}_{h_s}(c)}\sum_{\mbf{i}'\in \mb{S}_{h_s'}(c)}\sum_{t=1}^{r_3-h_t}\sum_{t'=1}^{r_3-h_t'}\frac{\PP\left(u_n^{-1}(X_t(c),X_{t'}(c))\in A^2 , (X_{t+h_t}(\mbf{i}),X_{t'+h_t'}(\mbf{i}'))\in B^2\right)}{rv_n \cdot \#\mb{S}_{h_s}(c)\cdot \#\mb{S}_{h_s'}(c)}
\end{equation*}
and 
\begin{equation*}
	\sum_{\mbf{i}\in \mb{S}_{h_s}(c)}\sum_{t=1}^{r_3-h_t}\sum_{t'=1}^{r_3}\frac{\PP\left(u_n^{-1}(X_t(c),X_{t'}(c))\in A^2 , X_{t+h_t}(\mbf{i})\in B\right)}{rv_n \cdot \#\mb{S}_{h_s}(c)}
\end{equation*} converge to $\sigma_{A,B}((h_s,h_t),(h_s',h_t'))$ and $\sigma'_{A,B}(h_s,h_t)$, respectively, where $r=r_1r_2r_3$ and $c=(\lceil (r_1+1)/2\rceil, \lceil (r_2+1)/2\rceil)$ (the ``center" of the block $B_{11}=[r_1]\times[r_2]$).
\\
\begin{itemize}
	\item[\textbf{(C)}]
	$\displaystyle\sum_{(c,t), (c',t')\in C(r_1,r_2)\times [n_3]} \PP\left(u_n^{-1}(X_t(c), X_{t'}(c')) \in \ A\times A \right)=\mc{O}(1)$,
\end{itemize}
where $C(r_1,r_2):=\{ c_{ij}\in [n_1]\times[n_2]: (i,j)\in [m_1]\times[m_2]\}$ is set of the ``centers" of the blocks $B_{ij}=[(i-1)r_1+1: ir_1]\times[(j-1)r_2+1: jr_2]$.
\begin{prop}[CLT for the iso-extremogram estimator]
	\label{conv_extremogram_c5}
	Assume that the following conditions hold for the  $\RR^k$-valued space-time process $X=\left\{ X_t(\mbf{s}) : (\mbf{s},t)\in \ZZ^2\times [0,\infty)\right\}$: 
	\begin{enumerate}
		\item The process $X$ is stationary in both space and time and satisfies the condition (I).
		\item The sequence $(u_n)$ is such that $(\ref{RV_c5})$ holds. Moreover, $r \ll v_n^{-1} \ll n$ and $\sqrt{nv_n}  \ll r  \ll nv_n r_3$, where $n=n_1n_2n_3$, $r=r_1r_2r_3$, $r_i  \ll n_i$ and $r_i=r_{n_i}\lfled \infty$, for $i=1,2,3$.
		\item The conditions (Cov') and (C) hold, and Lindeberg condition (Lin) is satisfied for the normalized variables $X_{\vn,(\mbf{s},t)}=u_n^{-1}X_t(\mbf{s})$ together with the family of cluster functionals $\mc{F}_{A,B}$. Moreover, for each $k\in\NN$, the coefficient $T^*_{\mbf{n},\mbf{t}}(\mbf{f}_k)$ defined in (\ref{T2_c5}) converges to zero as $\vn\to\infty$, for all $k-$tuple of cluster functionals $(f_1,\ldots, f_k)\in \mc{F}_{A,B}^k$ and all $\mbf{t}\in \RR^k$. The same assumption holds together with the family $\mc{F}_{A}:=\left\{ f_{A,A,0,0}\right\}$, which contains a single functional.
	\end{enumerate}
	Then, for each $(L_s,L_t)\in \NN_0\times \NN_0$, 
	\begin{equation}\label{TCL_extremogram_c5}
	\frac{\sqrt{n v_\vn}}{r_1r_2}\left( \h{\rho^*}_{A,B}(h_s, h_t) - \rho^*_{A,B,n}(h_s, h_t)\right)_{0\leq h_s\leq L_s, \ 0\leq h_t\leq L_t} \overset{\mc{D}}{\underset{\vn\to\infty}{\lfled}} \mc{N}(0,\Sigma_{A,B,L_s,L_t}),
	\end{equation}
	where $\rho^*_{A,B,n}(h_s, h_t):= \PP\left(u_n^{-1}X_{h_t}(h_s  \vec{e}_1)\in B \left| u_n^{-1}X_0(\mbf{0})\in A\right. \right)$ and $\Sigma_{A,B,L_s,L_t}$ is the covariance matrix, defined by the coefficients
	\begin{equation}\label{covarianza_ext_c5}
	\sigma_{\mbf{h},\mbf{h}'}=\sigma_{A,B}(\mbf{h},\mbf{h}')-\rho^*_{A,B}(\mbf{h}')\sigma'_{A,B}(\mbf{h})-\rho^*_{A,B}(\mbf{h})\sigma'_{A,B}(\mbf{h}')+\rho^*_{A,B}(\mbf{h})\rho^*_{A,B}(\mbf{h}')\sigma'_{A,A}(\mbf{0}),
	\nonumber
	\end{equation}
	with $\mbf{h},\mbf{h}'\in [0:L_s]\times[0:L_t]$.
\end{prop}
\textbf{Proof.}
	Consider the expression (\ref{est_extremogram_funct_c5}) of the iso-extremogram estimator. Then, for $(h_s,h_t)\in [0: L_s]\times [0,L_t]$, we obtain that
	%\vspace{-0.23cm}
	\begin{multline}\label{P5.1_1_c5}
	\frac{\sqrt{nv_n}}{r_1r_2} \left(\h{\rho^*}_{A,B}(h_s,h_t)-\rho^*_{A,B,n}(h_s,h_t) \right)
	\\
	=\dfrac{Z_\vn(f_{A,B,h_s, h_t})-\left(\frac{m h_t v_n}{\sqrt{nv_n}}+Z_\vn (f_{A,A,0,0})\right)\rho^*_{A,B,n}(h_s,h_t)+\frac{\delta_\vn}{\sqrt{nv_n}} + R}{\frac{r_1r_2}{\sqrt{nv_n}}Z_\vn(f_{A,A,0,0})+1+ \frac{r_1r_2 R_{A,A,0,0}}{n v_n}},
	\end{multline} 
	where $Z_\vn(\cdot)$ denotes the empirical process of cluster functionals (\ref{PE}). Besides, here $R=(nv_n)^{-1}\left( R_{A,B,h_s,h_t}-\rho^*_{A,B,n} R_{A,A,0,0}\right)$ and $m=m_1m_2m_3$.
	\\
	Now, notice that Chebyshev's inequality applied on the random variables $R$ and $r_1r_2 R_{A,A,0,0}/nv_n$ implies that these variables converge to zero in probability as $\vn\to\infty$. Similarly, applying Chebyshev's inequality together with the condition $\sqrt{nv_n}=o(r)$, we prove that $(nv_n)^{-1/2}\delta_\vn \overset{P}{\lfled} 0$, as $\mbf{n}\to\infty$. This last condition ($\sqrt{nv_n}=o(r)$) also guarantees that $mh_tv_n(nv_n)^{-1/2}\underset{\vn\to\infty}{\lfled}0$. Again, Chebyshev's inequality on the random variable $\frac{r_1r_2}{\sqrt{nv_n}}Z_\vn(f_{A,A,0,0})$, followed by the condition (C) and $r =o(nv_n r_3)$, implies that this converges to zero in probability as $\vn\to\infty$. 
	%\\
	Thus,  
	$$\frac{\sqrt{nv_n}}{r_1r_2} \left(\h{\rho^*}_{A,B}(h_s,h_t)-\rho^*_{A,B,n}(h_s,h_t) \right)=Z_\vn(f_{A,B,h_s, h_t})-\rho^*_{A,B,n}(h_s,h_t)Z_\vn (f_{A,A,0,0})+o(1).$$
	
	\
	
	The assumption 3 implies, from Theorem \ref{theo1_c5}, that $(Z_n(f_{A,B,h_s,h_t}))_{(h_s,h_t)\in [0:L_s]\times [0:L_t]}$ converges to a centered Gaussian r.v. with covariance matrix $(\sigma_{A,B}(\mbf{h},\mbf{h}'))_{\mbf{h},\mbf{h}'\in [0:L_s]\times [0:L_t]}$,  for each $(L_s,L_t)\in \NN^2_0$. Using the same argument, we prove that $Z_n(f_{A,A,0,0})$ converges to a centered Gaussian variable with variance $\sigma_{A,A}(\mbf{0},\mbf{0})$. 
	
	\
	
	Finally, considering the existence of $\sigma'_{A,B}$ in (Cov'), we obtain the result.\hspace{0.3cm}%$\blacksquare$ 
	\hfill$\square$\vskip2mm\hfill
%%%%%%%%%%%%%%%%%%%%%%%%%%%%%%%%%%%%%%%%%%%%%%%%%%%%%%%%%%%%%%%%%%%%%%%%%%%%%%%%%%%%%%%%%%%%%%%%%%%
%\bibliographystyle{plain}
%\bibliography{mybibfile}
\bibliographystyle{plain}

\end{document}